\documentclass[aip,jcp,reprint,showkeys,showpacs]{revtex4-1}

\usepackage{color}
\usepackage{epsfig}
\usepackage{graphicx}            

\usepackage[caption=false]{subfig}

\usepackage{array} 

\usepackage{amsmath}
\usepackage{latexsym}
\usepackage{amssymb}
\usepackage{amsfonts}

\usepackage{ulem}

\usepackage{pgfplots}
\usetikzlibrary{plotmarks}

\newcommand{\Frac}[2]%
     {\frac{\displaystyle #1}{\displaystyle #2}}

\newcommand{\cO}{\mathcal{O}}
\newcommand{\e}{\ensuremath{\mathrm{e}}}

\newcommand{\beq}{\begin{equation}}
\newcommand{\eeq}{\end{equation}}

\renewcommand{\Re}{\mbox{Re}}
\renewcommand{\Im}{\mbox{Im}}

\newcommand{\comment}[1]{\textcolor{blue}{(#1)}}
\newcommand{\del}[1]{\textcolor{red}{\sout{#1}}}
\newcommand{\structurecomment}[1]{}

\newlength\figurewidth 
\newlength\figureheight

\let\oldFootnote\footnote
\newcommand\nextToken\relax
\renewcommand\footnote[1]{%
    \oldFootnote{#1}\futurelet\nextToken\isFootnote}
\newcommand\isFootnote{%
    \ifx\footnote\nextToken\textsuperscript{,}\fi}

\DeclareMathOperator\sech{sech}

\begin{document}

\title{%
Solving the Schr\"odinger eigenvalue problem  by the
imaginary time propagation technique using splitting methods with complex
coefficients
}

%
%
%

%
%

\author{Philipp Bader}
\author{Sergio Blanes}
\affiliation{%
Instituto de Matem\'atica Multidisciplinar,
  Universitat Polit\`ecnica de Val\`encia, 46022 Valencia, Spain}
\author{Fernando Casas}
\affiliation{
Institut de Matem\`atiques i Aplicacions de Castell\'o, Universitat Jaume I, E12071 Castell\'on, Spain}


\begin{abstract}
The Schr\"odinger eigenvalue problem is solved with the imaginary time propagation technique. The separability of the Hamiltonian makes the problem suitable for the application of splitting methods.
High order fractional time steps of order greater than two necessarily have negative steps and can not be used for this class of diffusive problems.
However, there exist methods which use fractional  complex time steps with positive real parts which can be used with only a moderate increase in the computational cost.
We analyze the performance of this class of schemes and propose new methods which outperform the existing ones in most cases.
On the other hand, if the gradient of the potential is available, methods up to fourth order with real and positive coefficients exist.
We also explore this case and propose new methods as well as sixth-order methods with complex coefficients. In particular, highly optimized sixth-order schemes for near integrable systems using positive real part complex coefficients with and without modified potentials are presented.
A time-stepping variable order algorithm is proposed and numerical results show the enhanced efficiency of the new methods.

\vspace*{0.5cm}

\end{abstract}


\maketitle

\section{Introduction}
\structurecomment{Motivation of EV problem: whenever \textbf{positive coefficients} are needed}


\structurecomment{Define problem} We consider the eigenvalue
problem for the stationary Schr\"odinger equation (SE)
($\hbar=m=1$),

\beq \label{eq:Problem}
   H\phi_i(x) = E_i\phi_i(x), \qquad i=0,1,2,\ldots
\eeq
 where
 \beq   \label{eq:Problem2}
    H = T + V(x) = -\frac{1}{2}\Delta + V(x),
\eeq
$V(x)$ denotes the interaction potential and $\Delta$ is the Laplacian operator.
%
Since the Hamiltonian $H$ is a Hermitian operator, its eigenvalues
$E_i$ are real valued, and its corresponding real eigenfunctions
$\phi_i(x)$ form a basis of the underlying Hilbert space.
This particular problem has attracted great interest among theorists and practitioners \cite{%
chin09aoi
,janecek08afa
,roy01tdq
} due to its relevance for the understanding of the atomic and
molecular structure of matter.

A widely used approach to solve this problem is based on using the
corresponding time-dependent Schr\"odinger equation in imaginary
time ($t = -i \tau$), whose formal solution is given by the
evolution operator $\exp(-\tau H)$. In this way, in general, any
initial condition, under the action of $\exp(-\tau H)$, converges
asymptotically to the ground state solution when $\tau \rightarrow
\infty$. Notice that the evolution operator $\exp(-\tau H)$ has
the same eigenfunctions as the problem
(\ref{eq:Problem})-(\ref{eq:Problem2}). This technique is usually
referred to as the imaginary time propagation method (ITP for
short). In this setting, only the action of $\exp(-\tau H)$ on a
wave function has to be computed \cite{auer01afo,lehtovaara07sot}.

The ITP method can be regarded as an analog of the well-known power method in numerical linear algebra \cite{trefethen97nla}. In this sense,
one may also consider the inverse power method: instead of the iterative application of the exponential operator $\exp(-\tau H)$,
the scheme $v_{n+1}=(H-\tilde E_i)^{-1}v_n$, $n = 0,1,2,\ldots$ is used for some given $\tilde E_i$. This iteration is known to converge
after normalization to the eigenvector with eigenvalue closest to $\tilde E_i$.
Although faster convergence than for the ITP method can be observed for an accurate initial guess $\tilde E_i \approx E_i$, in general, the algorithm needs more iterations until convergence \cite{aichinger05afc}.

Since the operators $\e^{-\tau V}$ and  $\e^{-\tau T}$ can be exactly computed in the coordinate and momentum space, respectively,
the operator splitting technique
involving a composition of these exponential operators with appropriate coefficients can be used to approximate $\e^{-\tau H}$.
The computational cost depends on the number of changes between these coordinates which are cheaply performed by Fast Fourier transforms (FFT).

However, the operator splitting technique has some limitations. In particular,
splitting methods of order $p>2$ require negative time-steps \cite{sheng89slp,suzuki91gto}
and the instabilities caused thereof are analogous to the ones for the integration of a diffusion equation backwards in time.
 If it is feasible to compute the gradient of the potential $V$,
generalized splitting methods allow to build methods with positive
coefficients up to fourth order
\cite{chin97sif,omelyan02otc,chin02gsa}, but higher order methods
also use negative time-steps. In this paper we propose methods to
overcome the order barriers for both cases by using complex
time-steps. Splitting methods can be tailored to particular
equations to achieve better performances and we present criteria
based on near-integrability that apply to a wide range of problems
and thus yield highly efficient high order schemes. The obtained
methods outperform  the existing splitting schemes when high
accuracy is desired and could be appropriate for elaborating a
variable order algorithm. We also report some numerical
experiments illustrating the efficiency of the new methods.

\section{The imaginary time integration method for the Schr\"odinger equation}

\structurecomment{eigenvalue problem}%
An important property of the Hermitian operator $H$ is that (choosing properly the origin of the potential) its
eigenvalues $0\leq E_0\leq E_1\leq \dots $ are real and
nonnegative, and the corresponding eigenfunctions $\phi_i$ can be
chosen to form a real orthonormal basis on its domain.
The problem (\ref{eq:Problem}) originates from the time-dependent SE
 \beq \label{eq:SE}
    i\frac{\partial}{\partial t} \psi(x,t) = H\psi(x,t),
    \qquad \psi(x,0)=\psi_0(x).
 \eeq
A Wick rotation of the time coordinate, $t=-i\tau$, transforms \eqref{eq:SE} into a diffusion type equation
\beq\label{eq:SEimag}
    -\frac{\partial}{\partial \tau} \psi(x,\tau)
    = H\psi(x,\tau)
    , \qquad \psi(x,0)=\psi_0(x),
\eeq
with formal solution $\psi(x,\tau) = e^{-\tau H}\psi(x,0)$.
After expanding the initial condition $\psi_0$ in the basis of eigenfunctions $\phi_i$,
\[
  \psi_0(x) = \sum_i c_i\;\phi_i(x), \qquad
  c_i=\left\langle\phi_i(x)\,|\,\psi(x,0)\right\rangle,
\]
where $\langle\cdot\,|\,\cdot\rangle$ is the usual $L^2$ scalar
product, the time evolution of \eqref{eq:SEimag} is given by
 \beq \label{eq:evolution_imagtime}
    \psi(x,\tau) = e^{-\tau H}\psi(x,0) = \sum_i e^{-\tau
    E_i}\,c_i\;\phi_i(x).
\eeq
Asymptotically, for a sufficiently long time integration, we get $\psi(x,\tau)\to e^{-\tau E_0}\,c_0\phi_0$ since the other exponentials decay more rapidly. The convergence rate depends of course on the separation of the eigenvalues. For simplicity,
we restrict ourselves to the non-degenerate case $E_0<E_1$.
If there is degeneracy, it converges to a linear combination of eigenfunctions, and repeating this process with different initial conditions one can obtain a complete set of independent vectors of the subspace which can be orthonormalized.

Normalization of the asymptotic value yields the eigenfunction
$\phi_0$ and the corresponding eigenvalue is computed via $E_0=
\langle \phi_0|H\phi_0\rangle$. Excited states can be obtained by
propagating different wave functions simultaneously (or
successively) in time and using, for example, the Gram-Schmidt
orthonormalization or diagonalizing the overlap matrix
\cite{aichinger05afc}.

For simplicity in the presentation, the spatial dimension is set
to one unless it is explicitly stated, but our results also apply
to higher dimensions.

The problem is further simplified by assuming $x\in[a,b]$ with the interval $[a,b]$ sufficiently large such that the wave function and all its derivatives of interest vanish at the boundaries.
For numerical computations, the infinite dimensional domain of $H$ has to be truncated, which is done by discretizing the spatial coordinate $x$:
we fix $N$ equally spaced grid points $x_i=x_0+k\Delta x, \ k=0,1,2,\ldots,N-1$, with $a=x_0$ and $b=x_N$.
In this way, the interval is divided into $N$ subintervals of size $\Delta x=(b-a)/N$.

The potential $V$ is represented in this grid by a diagonal matrix
and the periodicity of the system ($\psi^{(n)}(a)=\psi^{(n)}(b)=0,
\ n=0,1,2,\ldots$) allows for the use of spectral methods (in
space) for the calculation of $T$, namely the Fast Fourier
Transform after which the matrix representation of $T$ also
becomes diagonal. The computational costs for the application of
$V$ and $T$ to a vector are thus proportional to $N$ and $N\log N$
operations, respectively. In a $d$-dimensional space with $N$ mesh
points on each dimension, their costs are proportional to $N^d$
and $N^d\log N$, respectively.


\section{Splitting methods for the Schr\"odinger equation}

To approximate the time evolution
\eqref{eq:evolution_imagtime}, i.e., the computation of $\e^{-\tau H}$
acting on a vector, we propose to use compositions of the operators $\e^{-\tau V}$ and
$\e^{-\tau T}$ evaluated at different times. A first example is provided by the well-known 
Strang splitting
\begin{equation}\label{eq:leapfrog}
  \Psi^{[2]}_h \equiv \e^{-\frac{h}2 V }\, \e^{-h T }\, \e^{-\frac{h}2 V },
\end{equation}
verifying $\Psi^{[2]}_h
=\e^{-hH}+\mathcal{O}(h^3)$ with $h\equiv\Delta \tau$.
Higher order approximations can be obtained by a more general
composition
\begin{equation}\label{eq:SplitGen}
  \Psi^{[p]}_h \equiv \prod_{i=1}^m \e^{-a_ih T }\, \e^{-b_ih V },
\end{equation}
where $\Psi^{[p]}_h=\e^{-hH}+\mathcal{O}(h^{p+1})$ if the
coefficients $a_i,b_i$ are chosen such that they satisfy a number
of order conditions (with $m$ sufficiently large).
\structurecomment{limitations to second order} It is well-known,
however, that methods of order greater than two ($p>2$)
necessarily have negative coefficients
\cite{sheng89slp,suzuki91gto,goldman96nos} (a simple proof can
be found in Ref.~\cite{blanes05otn}).
While this is usually not a problem for the coefficients $b_i$,
having negative $a_i$ coefficients makes the algorithm badly
conditioned (in the limit $N\rightarrow\infty$).

Composition methods with coefficients $b_i$ positive are also convenient for unbounded potentials, e.g., $V(x)=x^2$, since negative values of $b_i$ can generate large roundoff errors in the exponential $\e^{-b_iV}$ at the boundaries if the interval-size of the spatial discretization is not appropriately chosen and the potential takes exceedingly large values.

Splitting methods are particularly appropriate for the numerical
integration of this problem since the choice of the time step,
$h$, is not affected by the mesh size. Taking a finer mesh (i.e.,
a larger value of $N$) does not necessarily lead to a smaller time
step, and the extra computational effort originates only from the
FFTs, whose cost is $N\log(N)$ (or $N^d\log(N)$ in a
$d$-dimensional problem with $N$ points on each coordinate).

One possible approach to derive the order conditions to be satisfied by the coefficients $a_i$, $b_i$
consists in applying the Baker-Campbell-Hausdorff formula to the
composition (\ref{eq:SplitGen}), which we assume consistent
$(\sum_ia_i=\sum_ib_i=1)$ \cite{hairer06gni}. Thus we get $ \Psi^{[p]}_h = \exp(-h \mathcal{H})$, with
\begin{align}\label{eq:BCH1}
 \mathcal{H}  = & \nonumber
   T + V + h f_{2,1} [T,V] \\
  &+ h^2 \big(f_{3,1} [T,[T,V]] + f_{3,2} [V,[T,V]]\big) + \cdots,
\end{align}
where $f_{i,j}$ are polynomials of degree $i$ in the coefficients
$a_k,b_k$ and the symbol $[T,V]$ stands for the commutator of the operators $T$ and $V$.
 Condition $f_{2,1}=0$ leads to second order methods,
and this can always be achieved by taking a left-right symmetric composition in (\ref{eq:SplitGen})
because all even terms automatically vanish. Methods of higher orders
require in addition $f_{3,1}=f_{3,2}=0$. Taking into
account consistency, these equations can be written as \cite{sema08}
\begin{eqnarray}
 f_{3,1} & : & \sum_{1\leq i<j\leq k\leq m} a_ib_ja_k = \frac16,
 \label{eq:f31}
 \\
 f_{3,2} & : &  \sum_{1\leq i \leq j\leq k\leq m+1} b_ia_jb_k = \frac16.
 \label{eq:f32}
\end{eqnarray}
These two conditions imply that at least one of the $a_i$ as well
as one of the $b_i$ become negative (see \cite{blanes05otn} and references therein), so that only methods of order two can be used for this problem.

There are several possibilities to circumvent this limitation, and in the following, we enumerate some of them.

\paragraph{Modified potentials.} If the kinetic energy operator in \eqref{eq:SEimag} is quadratic in momenta, then the
nested commutator 
\begin{equation}\label{eq:VVT}
  [V,[T,V]] = \left(\nabla V\left(x\right)\right)^T\left(\nabla V\left(x\right)\right)
\end{equation}
is diagonal in coordinate space. For this reason, \eqref{eq:VVT} is usually called {modifying} potential.
In consequence, $[V,[V,[T,V]]] = 0$ and we can replace
 the terms $\e^{-b_ih V }$ in (\ref{eq:SplitGen})
by the more general operator 
$$
	\e^{-b_ih V  - c_ih^3[V,[T,V]]} 
$$
involving two parameters. As a result,
condition (\ref{eq:f32}) becomes
\begin{equation}\label{}
 f_{3,2} \ : \  \sum_{1\leq i \leq j\leq k\leq m+1} b_ia_jb_k
 + \sum_{i=1}^m c_i= \frac16.
\end{equation}
%

This equation can always be satisfied with a proper choice of the
coefficients $c_i$, so that the constraints on the coefficients
$a_i,b_i$ reduce to the single condition $f_{3,1}=0$, allowing for
positive coefficients. In addition, solutions with positive $c_i$
coefficients also exist. A first example is the 4th-order composition
\cite{koseleff94fcf,chin97sif}
\begin{equation}\label{eq:Chin4}
  \Psi^{[4]}_h \equiv
  \e^{-\frac{h}6 V }\, \e^{-\frac{h}2 T }\,
  \e^{-\frac{2h}3 V - \frac{h^3}{72} [V,[T,V]] }\,
  \e^{-\frac{h}2 T }\, \e^{-\frac{h}6 V }.
\end{equation}
It turns out, however, that 6th-order methods using the operator \eqref{eq:VVT} necessarily have
some negative coefficients $a_i$ \cite{chin05sop}.

\paragraph{Near-integrable systems.} When the Hamiltonian can be considered as a perturbed system,
i.e.,  $H=H_0+\varepsilon V_{\varepsilon}(x)$  with an exactly
solvable part $H_0=T+V_0(x)$ and a small perturbation $\varepsilon
V_{\varepsilon}(x)$, it is advantageous to split the Hamiltonian
into the dominant part $H_0$ and its perturbation $\varepsilon
V_{\varepsilon}$. For example, if one is interested in the lower
excited states, which evolve near the minimum of the potential, it
can be useful to separate the quadratic part and to treat the
remainder as a perturbation since the harmonic oscillator has a
simple and fast solution  using FFTs \cite{chin05foa,bader11fmf}.

Notice that in this case, the commutator 
\begin{equation*}
  [\varepsilon V_{\varepsilon},[H_0,\varepsilon V_{\varepsilon}]] =
  \varepsilon^2 \left(\nabla V_{\varepsilon}\left(x\right)\right)^T\left(\nabla
  V_{\varepsilon}\left(x\right)\right)
\end{equation*}
depends only on the coordinates and modified potentials can also be applied as before. Then, all
compositions remain the same except for  replacing $T$ by $H_0$
and $V$ by $\varepsilon V_{\varepsilon}$.


With the additional information that one part of the operator is significantly
smaller than the other, it is clear that the error expansion for a
consistent method $\Psi_h$ can be asymptotically expressed as
$$
    \Psi_h- \e^{-hH} = \sum_{i \ge 1}\sum_{k \ge s_i} e_{i,k}\,\varepsilon^i h^{k+1}, \text{ as } (h, \varepsilon) \rightarrow (0,0),
$$
where the $s_i$ start from the first non-vanishing error
coefficient $e_{s_i,k}$. We say that $\Psi_h$ is of generalized order
$(s_1,s_2,\ldots, s_m)$ (where $s_1 \ge s_2 \ge \cdots \ge s_m$) if the local error satisfies that
\[
   \Psi_h- \e^{-hH}  = \mathcal{O}(\varepsilon h^{s_1 + 1} + \varepsilon^2 h^{s_2 + 1} + \cdots + \varepsilon^m h^{s_m + 1}).
\]
Thus, for a method of generalized order $(8,2)$, denoted by $\Psi^{(8,2)}_h$, the error reads
\[
    \Psi^{(8,2)}_h - \e^{-hH} = e_{1,9}\varepsilon h^9 + e_{2,2}\varepsilon^2 h^3 + \mathcal{O}\left(\varepsilon^3 h^3\right).
\]
This class of schemes can also be applied in several other
situations. For instance, suppose one takes a sufficiently fine
mesh. Then $\|T\|\gg \|V\|$ and the previous considerations apply
(with $H_0 = T$). 
Also, if $V(x)=x^n$, then the virial theorem $\langle \phi\,|T|\,\phi\rangle = \langle \phi\,|\nabla V(x) x|\,\phi\rangle$ leads to $\langle
T\rangle=n\langle V\rangle$.

\paragraph{Complex coefficients.} A third possibility consists of considering complex coefficients in the composition  (\ref{eq:SplitGen})
(with or without modified potentials). In other problems where the presence
 of negative real coefficients is unacceptable,
the use of high-order splitting methods with complex coefficients having positive real part has shown to possess some advantages.
In recent years a systematic
search for new methods  with complex coefficients has been carried out and the resulting schemes have been tested in different settings:
Hamiltonian systems in celestial mechanics  \cite{chambers03siw},
the time-dependent Schr\"odinger equation in quantum mechanics \cite{bandrauk06cis}
and also in the more abstract setting of evolution equations with unbounded operators generating analytic semigroups \cite{castella09smw,hansen09hos}.
{It is worth noticing that the
propagator $\exp(z \Delta)$ ($z \in \mathbb{C}$) associated with the Laplacian
 is well-defined (in a reasonable distributional sense) if and only if $\Re(z) \ge 0$
 \cite{castella09smw}, which is the case for the presented methods.}


Many of the existing splitting methods with complex coefficients have been constructed by applying the composition technique to the symmetric second-order leapfrog scheme \eqref{eq:leapfrog}.
For example, a fourth-order integrator can be obtained with the symmetric composition
\begin{equation}\label{eq:TripleJump}
  \Psi_h^{[4]} = \Psi_{\alpha h}^{[2]} \ \Psi_{\beta h}^{[2]} \ \Psi_{\alpha h}^{[2]},
\end{equation}
where
\begin{equation}  \label{15pr}
 \alpha = \frac{1}{2-2^{1/3} \e^{2ik\pi/3}}
 , \qquad
 \beta = \frac{2^{1/3} \e^{2ik\pi/3}}{2-2^{1/3} \e^{2ik\pi/3}},
\end{equation}
and $k=1,2$. In both cases, one has $\Re(\alpha), \Re(\beta) >0$.
Higher order composition methods with complex coefficients and positive real part can be found in Refs.~\onlinecite{castella09smw,hansen09hos,blanes12oho}, where several numerical examples are also reported.


\section{New splitting methods for the ITP problem}

In this section, we carry out a systematic search of methods within the classes (a)-(c) above
enumerated. 
The best methods for each subclass can be found online\cite{website} with 25 digits of accuracy whereas the methods used in the numerical examples (Sec.~\ref{sec:numerics}) are given in the subsequent tables with 18 digits for simplicity.

We only consider symmetric methods and, since $T$ and $V$ have qualitatively different properties, we analyze both TVT-and VTV-type compositions, defined as
\begin{align}
 \Psi_h^{[p]} & =
                             \e^{-a_1hT} \e^{-b_1hV} \e^{-a_2hT} \cdots
 \e^{-a_2hT} \e^{-b_1hV} \e^{-a_1hT},    \label{eq:TVT} \\
 \text{ and } \nonumber \\
 \Psi_h^{[p]} & =
                             \e^{-b_1hV} \e^{-a_1hT} \e^{-b_2hV} \cdots
 \e^{-b_2hV} \e^{-a_1hT} \e^{-b_1hV}  \label{eq:VTV},
\end{align}
respectively.
In principle, both compositions have the same computational cost for the same number of exponentials. 
Nevertheless, due to a projection step to the real part after each full time-step, only in the VTV composition we can concatenate the last map in the current step with the first stage in the next one. 
The TVT compositions thus require two additional FFTs in comparison with the VTV composition, and this is accounted for in the numerical experiments.

The methods we obtain are classified into two families: (I) methods without modified potentials and (II) methods with modified potentials. 
For each class we distinguish between methods for general problems (with the unique constraint that $[V,[V,[T,V]]]=0$) and methods for near-integrable problems (when the main dominant part contains the kinetic energy).

We have explored both TVT and VTV compositions with different
number of stages. In some cases we consider extra stages to have
free parameters for optimization. When the number
and complexity of the order conditions is relatively low, we
get all solutions. We then select the solutions having
all of their coefficients with positive real part. Finally, we choose
the solution which minimizes
\begin{equation}\label{eq:SumCoefs}
   \sum_i ( |a_i| + |b_i|)
\end{equation}
and/or minimizes the absolute value of the real part of the
coefficients appearing at the leading error terms. These
methods are subsequently tested on several numerical examples. After this process, we
collect a number of schemes offering the best performance for most of the problems
considered. In practice, however, one has to bear in mind that the relative performance
between different methods depends eventually on the particular problem considered, the desired accuracy,
the initial conditions, etc.

\subsection{ Methods without modified potentials}

TVT and VTV compositions with 3 up to 9 stages have been analyzed. To simplify the notation, we denote compositions
\eqref{eq:TVT} and \eqref{eq:VTV} as
\begin{align*}
  {\rm T}n_m &= a_1 \, b_1 \, a_2  \, \cdots \, a_2 \, b_1 \, a_1
   ,\\ 
  {\rm V}n_m &= b_1 \,a_1 \,  b_2  \, \cdots \, b_2 \, a_1 \, b_1
\end{align*}
respectively. Here $n$ indicates the order (or generalized order) of the method and
$m$ corresponds to the number of stages, i.e., the number of $b_i$
coefficients in the TVT composition or the number of $a_i$
coefficients in the VTV composition.
The coefficients of the selected TVT methods are collected in Table~\ref{tab.1b}, whereas those corresponding to the TVT methods are displayed in Table~\ref{tab.1a}.

\subsubsection{Methods for general problems}
Analogously to \eqref{eq:BCH1}, the symmetric compositions \eqref{eq:TVT} and \eqref{eq:VTV} can be formally expressed as a single exponential $\Psi_h^{[p]} = \exp(-h \mathcal{H})$ with polynomials $f_{i,j}$ in $a_k,b_l$ multiplying commutators $E_{i,j}$:
\begin{align*}
 \mathcal{H} & =  {  T + V + h^2 \big(f_{3,1} E_{3,1} + f_{3,2} E_{3,2}\big)  }  \\
  &  {   + h^4 \big(f_{5,1} E_{5,1} + f_{5,2} E_{5,2} + f_{5,3} E_{5,3} + f_{5,4} E_{5,4}\big)  } \\
  &  {   + h^6 \big(f_{7,1} E_{7,1} + f_{7,2} E_{7,2} + \cdots \big) + \cdots, }
\end{align*}
where the $E_{i,j}$ are chosen to form a basis of the algebra 
of commutators of length $i$. 
The chosen basis elements relevant for our exposition are
\begin{align*}
    E_{3,1} &= [T,[T,V]], &
    E_{3,2} &= [V,[T,V]],\\
    E_{5,1} &= [T,[T,[T,[T,V]]]],&
    E_{5,2} &= [V,[T,[T,[T,V]]]],\\
    E_{5,3} &=-[T,[V,[T,[T,V]]]],&
    E_{5,4} &= [V,[V,[T,[T,V]]]],\\
    E_{7,1} &= [T,[T,E_{5,1}]],&
    E_{7,2} &= [V,[T,E_{5,1}]].
\end{align*}
Here we summarize some of the methods which have been analyzed:

\paragraph{\textit{3-stage compositions}.} A 3-stage composition has
 sufficient parameters to build 4th-order methods. There is one real
  solution and two complex solutions (conjugate to each other). For example, the
  VTV method corresponds to the composition (\ref{eq:TripleJump}) when
  $\Psi_{h}^{[2]}$ is given by (\ref{eq:leapfrog}). The TVT
  version is obtained by interchanging $T$ and $V$.

\paragraph{\textit{5-stage compositions}.} Fourth-order methods with
  two free parameters can be obtained using 5-stage symmetric
  compositions. These two parameters can be used to build methods
  of effective order 6 (i.e., 4th-order methods that are conjugate to
  6th-order methods by a near-identity change of variables). This requires to impose some additional constraints on
  the leading error terms, $f_{5,j}, \ j=1,2,3,4$. Specifically, these are
  $f_{5,1}-f_{5,2}=0$ and
  $f_{5,3}+f_{5,4}=0$ \cite{celesmech2000}.  We have found six solutions for the TVT composition and three solutions for the VTV composition with coefficients having positive real part.
The solutions with smallest error terms at order 5 are denoted by T4$_{5}$ and V4$_{5}$\cite{website}.

\paragraph{\textit{7-stage compositions}.} In principle, there are sufficient
parameters to build 6th-order methods with 7 stages. For the TVT composition there are 11 solutions with all coefficients having positive real parts. 
The solution leading to a minimum value of the norm of the error at order 7 can be found online\cite{website}.

With respect to the VTV composition, the best method we have found is identical with the most efficient sixth-order method obtained by Chambers \cite{chambers03siw}, where it has been presented as a symmetric composition similar to \eqref{eq:TripleJump} but with 7 stages instead of 3, and with $\Psi_{h}^{[2]}$  given by \eqref{eq:leapfrog}.



\subsubsection{Methods for near-integrable problems}

Proceeding analogously as before, we arrive at the following methods. We recall that in all compositions one should replace $T$ by $H_0$ and $V$ by $\varepsilon V_{\varepsilon}$.


\paragraph{\textit{$n$-stage compositions of generalized order $(2n,2)$}.}
This class of compositions has real and positive coefficients
\cite{McL95a,laskar01hos}. A 4-stage VTV composition of
generalized order $(8,2)$ is given by scheme
V84M$_{4}^{\text{LR}}$
 in Table~\ref{tab.2a} with $c_1=0$.

\paragraph{\textit{5-stage compositions}.} To build a method of generalized order
(8,4) the following conditions must be satisfied by a consistent
and symmetric method: $f_{3,1}=f_{3,2}=f_{5,1}=f_{7,1}=0$. It
requires at least 5 stages, and in this case only one solution
with all coefficients having positive real part is found both for
the TVT and VTV compositions. The coefficients of these
methods, denoted by T84$_{5}$ and V84$_{5}$, are collected in
Table~\ref{tab.1b} and Table~\ref{tab.1a}, respectively.

\paragraph{(8,6,4) \textit{methods}.} To build a  (8,6,4) method, the
coefficients of a consistent and symmetric method must satisfy the
following order conditions:
$f_{3,1}=f_{3,2}=f_{5,1}=f_{5,2}=f_{5,3}=f_{7,1}=0$. They therefore require
at least 7 stages. In this case, it is possible to get all solutions.
Scheme T864$_{7}$ corresponds to the solution minimizing (\ref{eq:SumCoefs}), whereas
V864$_{7}$ provides the minimum value of $|f_{5,3}+f_{5,4}|$.

\paragraph{(8,6) \textit{methods}.} Increasing the number of stages to 9 we have two free parameters, which are used to
satisfy in addition the following conditions: $f_{5,4}=f_{7,2}=0$.  In this way, methods of generalized
order (8,6) and effective order (10,8,6) are obtained.
Two efficient schemes 
correspond to T86$_{9}$ and
V86$_{9}$ in Table~\ref{tab.1b} and Table~\ref{tab.1a}, respectively \cite{Ander_Joseba}.

\subsection{Methods with modified potentials}

Fourth-order methods incorporating modified potentials do exist
with real and positive coefficients. In fact, 2- and 3-stage
schemes have been extensively studied\cite{chin05sop,omelyan02otc,chin02gsa}. 
Methods of generalized order $(n,4)$ also exist with positive real coefficients\cite{laskar01hos}. 
Here we construct new methods of generalized order (6,4) and (8,4) with this property and generalize the treatment to 6th-order schemes with complex coefficients. In all
cases, we take compositions TVT and VTV with up to 5 stages and denote them as
\begin{eqnarray}\label{}
  {\rm T}n{\rm M}_m &=&  a_1 \, (b_1 \, c_1) \, a_2   \, \cdots \, a_2 \,  (b_1 \, c_1) \, a_1,
 \nonumber \\
  {\rm V}n{\rm M}_m &=&  (b_1 \, c_1) \, a_1 \, (b_2 \, c_2)  \, \cdots \,
  (b_2 \, c_2) \, a_1  \, (b_1 \, c_1).  \nonumber
\end{eqnarray}
Here, the parenthesis is used to help counting of the number of
exponentials, and the letter M indicates that the methods use
modified potentials. Notice that the number of free parameters can
differ for the TVT and VTV sequences with the same number of
exponentials because the exponent of a modified potential contains
two parameters. The coefficients of the selected methods are
collected in Table~\ref{tab.2b} and  Table~\ref{tab.2a} for the
TVT and VTV compositions, respectively.

\subsubsection{Methods for general problems}


\paragraph{\textit{4-stage compositions}.}
Under the restriction of having real positive coefficients, we have obtained the fourth-order VTV method OMF-4M, already discovered in Ref.~\onlinecite{omelyan02otc} (eq. (36) therein).

The VTV composition allows one to build 6th-order methods, whereas the TVT needs an extra stage. 
There is only one solution (and its complex conjugate) with all coefficients having positive real part.
It is denoted by V6M$_{4}$ and can be found online\cite{website}.


\subsubsection{Methods for near-integrable problems.}

We first consider $(n,4)$ methods with real
and positive coefficients. For schemes of generalized order (8,6) we collect only complex solutions
with positive real part.

\Roman{paragraph}
\paragraph{(6,4) \textit{methods}}
They require at least 3 stages to satisfy the following order conditions: $f_{3,1}=f_{3,2}=f_{5,1}=0$. 
The coefficients $a_i$ and $b_i$ correspond to the methods (6,2) obtained in Ref.~\onlinecite{McL95a} (without modified potentials).
We have also considered methods with 4 stages in order to have additional free parameters.
As previously mentioned, there is the same number of order conditions as parameters to get a method of order 6 for the VTV sequence, but there are no solutions with coefficients being real and positive. 
To get a sixth-order method the following conditions must also to be satisfied: $f_{5,2}=f_{5,3}=f_{5,4}=0$. 
The coefficients $c_i$ only appear in $f_{5,3}$ and $f_{5,4}$ and can only be used to cancel these terms. 
The VTV sequence has three free parameters which can be used to annihilate $f_{5,3}$ and $f_{5,4}$ and to minimize the absolute value of $f_{5,2}$ under the constraint that all coefficients must be real and positive. 
The TVT sequence has only two free parameters which can be used to annihilate $f_{5,3}$ and to minimize the absolute value of the dominant term, $f_{5,2}$, under the same constraint on the coefficients.
The best methods we have obtained are denoted by T64M$_{4}$ and V64M$_{4}$ and are published online\cite{website}.

\paragraph{(8,4) \textit{methods}}
They require at least 4 stages. The coefficients
$a_i$ and $b_i$ correspond to the methods (8,2) without using
modified potentials and obtained in \cite{McL95a}. There is one
coefficients $c_i$ in the TVT composition  which can be used to
cancel $f_{5,3}$, and two coefficients $c_i$ in the VTV
composition which can be used to annihilate $f_{5,3}$ and $f_{5,4}$.
The solution with $c_2=c_3=0$ was already obtained in
\cite{laskar01hos}. We have collected the corresponding
coefficients for this method, V84M$_{4}^{\text{LR}}$, in
Table~\ref{tab.2a}. We have also considered methods with 5 stages
in order to have an additional free parameters. There is the same
number of order conditions as parameters to get a method of order
(8,6) (which would be of order 6 for a general problem) but,
obviously, there are no solutions with coefficients real and
positive.
As in the previous case, the term  $f_{5,2}$ can not be
zeroed using real positive coefficients. Then  in both TVT and
VTV compositions we have chosen the method which, while having
real and positive coefficients, minimize its absolute value. The
best methods we have obtained are denoted by T84M$_{5}$ and
V84M$_{5}$.

\paragraph{(8,6) \textit{methods}}
They require at least 5 stages and do not
admit real and positive solutions for the coefficients and we are
forced to consider complex solutions. We have found only one
solution with positive real part in the coefficients for both TVT
and VTV compositions. The coefficients for the methods denoted
by T86M$_{5}$ and V86M$_{5}$ are given in Table~\ref{tab.2b} and
Table~\ref{tab.2a}, respectively.


\begin{table}[tbp]
\caption{Compositions TVT without modified potentials.}
\label{tab.1b} {
\begin{tabular}{l}
\begin{tabular}{l}
\hline
    \\
T84$_{5}= a_1 \, b_1 \, a_2  \,b_2 \,
 { a_3 \, b_3 \, a_3} \, b_2 \, a_2 \, b_1 \, a_1 $
 \\
\hline\hline
   $a_1=0.071401131540044698 + 0.010155431019886789 i$ \\
   $b_1=0.178696854264631978 + 0.028197506313218021 i$ \\
   $a_2=0.236383805190074736 + 0.070427007139534522 i$ \\
   $b_2=0.198453474708154649 + 0.082962314733854963 i$ \\
   $a_3=1/2-(a_1+a_2) = 0.1922... - 0.0806...i$  \\
   $b_3=1-2(b_1+b_2)  = 0.2457... - 0.2223...i$ \\
\end{tabular}
\\
\begin{tabular}{l}
\hline
    \\
T864$_{7}=
 a_1 \, b_1 \, a_2 \, b_2 \, a_3 \, b_3 \,
  { a_4 \, b_4 \, a_4}
 \,  b_3 \, a_3 \, b_2 \, a_2 \, b_1 \, a_1$
  \\
\hline\hline
   $a_1=0.055705821110864236 + 0.018670384565085049 i $ \\
   $b_1=0.115779449626990422 + 0.046131356173382847 i$ \\
   $a_2=0.118843282163492564 -  0.024151805322796634 i$ \\
   $b_2=0.129128920804026450 -  0.119039413303774209 i$ \\
   $a_3=0.158591515575195578 -  0.076302551893579599 i$ \\
   $b_3=0.184643464154438944 -  0.003053761445376182 i$ \\
   $a_4 = 1/2 - (a_1+a_2+a_3) =  0.1669... + 0.0818...i $\\
   $b_4 = 1 - 2(b_1+b_2+b_3)  = 0.1409... + 0.1519...i$ \\
   \hline
\end{tabular}
\\
\begin{tabular}{l}
\hline
    \\
T86$_{9}=
 a_1 \, b_1 \, a_2 \, b_2 \, a_3 \, b_3 \, a_4 \, b_4 \,
  { a_5 \, b_5 \, a_5}
 \, b_4 \, a_4 \, b_3 \, a_3 \, b_2 \, a_2 \, b_1 \, a_1$
  \\
\hline\hline
   $a_1=0.042257897299860339 - 0.014215780224181831 i$ \\
   $b_1=0.094894869367770736 - 0.037963806472588094 i$ \\
   $a_2=0.095260398471830494 + 0.004518725891475591 i$ \\
   $b_2=0.097374660381711248 + 0.088518877931710497 i$ \\
   $a_3=0.099960578944766657 + 0.090271995071312563 i$ \\
   $b_3=0.118584793520055816 + 0.038356250608401259 i$ \\
   $a_4=0.148695530402608487 + 0.011438117187614089 i$\\
   $b_4=0.136865119760326031 - 0.023587404969570006 i$ \\
   $a_5=1/2-(a_1+a_2+a_3+a_4) = 0.1138... - 0.0920...i$  \\
   $b_5=1-2(b_1+b_2+b_3+b_4)  = 0.1046... - 0.1306...i$ \\
   \hline
\end{tabular}
\end{tabular}
}
\end{table}

\begin{table}[tbp]
\caption{Compositions VTV without modified potentials.}
\label{tab.1a} {
\begin{tabular}{l}
\begin{tabular}{l}
\hline
    \\
V84$_{5}= b_1 \, a_1 \, b_2 \, a_2  \,
 { b_3 \, a_3 \, b_3} \, a_2 \, b_2 \, a_1 \, b_1 $
 \\
\hline\hline
   $b_1=0.052472525516129026 - 0.010958940842458138 i $\\
   $a_1=0.175962140656732362 - 0.054483056228160557 i$ \\
   $b_2=0.246023563332753880 - 0.125228547924834352 i$ \\
   $a_2=0.181259898687454283 - 0.034864508232090522 i$ \\
   $b_3=1/2-(b_1+b_2)= 0.2015... + 0.1362...i$  \\
   $a_3=1-2(a_1+a_2) = 0.2856... + 0.1787...i$ \\
\end{tabular}
\\
\begin{tabular}{l}
\hline
    \\
V864$_{7}=
 b_1 \, a_1 \, b_2 \, a_2 \, b_3 \, a_3 \,
  { b_4 \, a_4 \, b_4}
 \,  a_3 \, b_3 \, a_2 \, b_2 \, a_1 \, b_1$
  \\
\hline\hline
   $b_1=0.060017770752528926 -  0.009696150746907738 i$ \\
   $a_1=0.108904710931114447 -  0.075700232434276860 i$ \\
   $b_2=0.067017987316853817 + 0.003927567742822542 i$ \\
   $a_2=0.106594114300156182 + 0.139651903644940761 i$ \\
   $b_3=0.189300872388005476 + 0.091055103879530385 i$ \\
   $a_3=0.204897016414416105 + 0.009719057955143112 i$ \\
   $b_4=1/2-(b_1+b_2+b_3) = 0.1837... - 0.0853...i$\\
   $a_4=1-2(a_1+a_2+a_3)  = 0.1592... - 0.1473...i$\\
   \hline
\end{tabular}
\\
\begin{tabular}{l}
\hline
    \\
V86$_{9}=
 b_1 \, a_1 \, b_2 \, a_2 \, b_3 \, a_3 \, b_4 \, a_4 \,
  { b_5 \, a_5 \, b_5}
 \, a_4 \, b_4 \, a_3 \, b_3 \, a_2 \, b_2 \, a_1 \, b_1$
  \\
\hline\hline
   $b_1=0.032497706037458608 + 0.010641310380458924 i$ \\
   $a_1=0.087895680441261752 + 0.036052576182866484 i $ \\
   $b_2=0.094180923422602148 + 0.023866875362648754 i$ \\
   $a_2=0.095351855399045611 - 0.065128376035135147 i$ \\
   $b_3=0.101132953097231180 - 0.112201757337044841 i$ \\
   $a_3=0.121865575594908413 - 0.054974002471495827 i$ \\
   $b_4=0.160941382119434892 - 0.016127643896952891 i$ \\
   $a_4=0.141506882718462097 + 0.024607229046524026 i$\\
   $b_5=1/2-(b_1+b_2+b_3+b_4) = 0.1112... + 0.0938...i$  \\
   $a_5=1-2(a_1+a_2+a_3+a_4)  = 0.1068... + 0.1189...i$ \\
   \hline
\end{tabular}
\end{tabular}
}
\end{table}

\begin{table}[tbp]
\caption{Compositions TVT with modified potentials.} \label{tab.2b} {
\begin{tabular}{l}
\begin{tabular}{l}
\hline
    \\
T84M$_{5}=   a_1 \, (b_1 \, c_1) \, a_2 \,(b_2 \, c_2)
  \, { a_3}  \, ({ b_3} \, c_3)  \,  { a_3} \, (b_2 \, c_2) \, a_2 \, (b_1 \, c_1) \, a_1$
\\
\hline\hline
   $a_1=0.058520963359694865$ \\
   $b_1=0.145381537601615725, \qquad  c_1=0.000245906549261228$ \\
   $a_2=0.207903047442871771$ \\
   $b_2=0.244351408696638327 , \qquad c_2=0.000259178561419125$ \\
   $a_3=1/2-(a_1+a_2)=0.2336...$  \\
   $b_3=1-2(b_1+b_2)= 0.2205...,   \quad  c_3=0.000938105701711153$ \\
\end{tabular}
\\
\begin{tabular}{l}
\hline
    \\
T86M$_{5}=    a_1 \, (b_1 \, c_1) \, a_2 \,(b_2 \, c_2)
  \, { a_3}  \, ({ b_3} \, c_3)  \,  { a_3} \, (b_2 \, c_2) \, a_2 \, (b_1 \, c_1) \, a_1$
\\
\hline\hline
   $a_1=0.063556051997493102 + 0.010606890396680920 i$ \\
   $b_1=0.156939525347224563 + 0.027931306200415819 i$ \\
   $c_1=0.000133739181746125 + 0.000085540153220213 i$ \\
   $a_2=0.208998817231756322 + 0.040240203826523395 i$ \\
   $b_2=0.222383136675982213 + 0.026033262090035938 i$ \\
   $c_2=0.000484323504408882 + 0.000241671051573332 i$ \\
   $a_3=1/2-(a_1+a_2) = 0.2274... - 0.0508...i$  \\
   $b_3=1-2(b_1+b_2)  = 0.2414... - 0.1079...i$ \\
   $c_3=0.000179180363327321 - 0.000858304413034511 i$ \\
   \hline
\end{tabular}
\end{tabular}
}
\end{table}

\begin{table}[tbp]
\caption{Compositions VTV with modified potentials.} \label{tab.2a} {
\begin{tabular}{l}
\begin{tabular}{l}
\hline
    \\
V84M$_{5}=  (b_1 \, c_1) \, a_1 \, (b_2 \, c_2) \, a_2  \,
 ({ b_3} \, c_3) \, { a_3} ({ b_3} \, c_3) \, a_2 \, (b_2 \, c_2) \, a_1 \, (b_1 \, c_1)$
\\
\hline\hline
   $b_1=0.042308451243127365, \qquad c_1=0.000232966269565498$ \\
   $a_1=0.142939324267716184$ \\
   $b_2=0.219303568753387110 , \qquad c_2=5.56677120231130 \cdot 10^{-7}$ \\
   $a_2=0.242474508234531493$   \\
   $b_3=1/2-(b_1+b_2)=0.2292..., \quad c_3=0.000794490777479431$  \\
   $a_3=1-2(a_1+a_2) =0.2384...$ \\
\end{tabular}
\\
\begin{tabular}{l}
\hline
    \\
V84M$_{4}^{\text{LR}}=  (b_1 \, c_1) \, a_1 \, (b_2 \, c_2) \, {\bf a_2
\,
 (b_3} \, c_3) \, {\bf a_2} \, (b_2 \, c_2) \, a_1 \, (b_1 \, c_1)$
\\
\hline\hline
   $b_1=1/20, \qquad  c_1=\frac{3861 - 791 \sqrt{21}}{129600}, \qquad a_1=1/2-\sqrt{3/28}$ \\
   $b_2=49/180 , \qquad c_2=0$  \\
   $a_2=1/2-a_1=\sqrt{3/28}$  \\
   $b_3=1-2(b_1+b_2)=16/45, \quad c_3=0$ \\
\end{tabular}
\\
\begin{tabular}{l}
\hline
    \\
V86M$_{5}=  (b_1 \, c_1) \, a_1 \, (b_2 \, c_2) \, a_2  \,
 ({\bf b_3} \, c_3) \, {\bf a_3} ({\bf b_3} \, c_3) \, a_2 \, (b_2 \, c_2) \, a_1 \, (b_1 \, c_1)$
\\
\hline\hline
   $b_1=0.046213625838152095 - 0.007824529355983108 i$ \\
   $c_1=0.000035830461339520 + 0.000074370857685421 i$ \\
   $a_1=0.152650950104799817 - 0.030279967163699065 i$ \\
   $b_2=0.224258052678856384 - 0.050879282402761772 i$ \\
   $c_2=0.000338053435041382 - 0.000490508913279372 i$ \\
   $a_2=0.226364275186039762 - 0.016537249619936515 i$ \\
   $b_3=1/2-(b_1+b_2)=  0.2295... + 0.0587...i$  \\
   $c_3=0.000408311644874003 + 0.000484371967433683 i$ \\
   $a_3=1-2(a_1+a_2)=0.2420... + 0.0936...i$ \\
\end{tabular}
\end{tabular}
}
\end{table}

\section{Numerical examples}\label{sec:numerics}
\subsection{Efficiency of the methods}
\structurecomment{The potentials} As test bench for the numerical
methods, we consider in the following two qualitatively different
cases, the P\"oschl-Teller potential and a perturbed harmonic
oscillator, the latter being a classic example of a
near-integrable system and of practical interest \cite{roy01tdq}.
These two problems can be numerically integrated using modified
potentials. However, we compare the relative performance of the
methods (with and without modified potentials) separately in order
to study the performance of the methods when it is not feasible to
compute the gradient of the potential.

\structurecomment{The procedure} The numerical integration
proceeds as follows: starting from random initial data, we iterate
with fixed time-step until the sufficiently large final time
$T=100$ and compare the result with the exact solution, $\psi(T)$,
which has been obtained by integrating with a much smaller time
step. The spatial interval is fixed for all experiments to
$[-10,10]$ and is discretized with 
 {$N=128$} equidistant mesh points.
Similar results are obtained for larger $N=256,512,1024$.
At each step, we project the obtained vector to its real part and
normalize it to one in $\ell_2(\mathbb{R})$,
i.e., given the method $\Psi_h^{[p]}$ and initial conditions,
$u_n\in\mathbb{R}^N$, we compute $u_{n+1}$ as
\[
  \tilde u_{n+1}=\Psi_h^{[p]} \, u_{n};
\]
then, since $\tilde u_{n+1}$ is a complex vector (but
$\mathcal{O}(h^p)$ away from a real vector) we project on the real
space by removing the imaginary part
\[
  \bar u_{n+1}= \Re(\tilde u_{n+1})
\]
and then normalize the solution
$
  u_{n+1}= {\bar u_{n+1}}\big/{\|\bar u_{n+1}\|},
$
where the norm is given by
$$
   \|w\|^2 \equiv \Delta x \sum_{j=0}^{N-1} w_j^2,\quad w=(w_0,\ldots,w_{N-1})\in \mathbb{R}^N.
$$
We take as the computational cost the number of Fourier transforms
necessary until the final time.
%
%
In addition, the methods using complex coefficients are penalized
by a factor $2$ in the computational cost, which comes from the
use of complex Fourier transforms instead of real FFT.
We repeat the numerical integrations for different values of the
time step, i.e., $h=T/M$ for different values of $M$. We take as
the approximate solution, $\phi(T)=u_n$ in each case and measure
the error as
$$
    \text{error} =\|\psi(T)-\phi(T)\|.
$$

\structurecomment{why we do this} This procedure will allow us to
determine the efficiency of the new splitting methods, which will
depend on the desired accuracy, and thereby choose the methods which
are most appropriate for implementation with a more efficient
algorithm that is based on variable time step and order.
\structurecomment{two types of methods} 
We distinguish two types of problems: on the one hand, methods that include modified potentials, the reference methods being Chin-4M \eqref{eq:Chin4}, OMF-4M \cite{omelyan02otc} and V84M$_4^{\text{LR}}$ \cite{laskar01hos} given in Table~\ref{tab.2a} as well as a differently optimized scheme SCF-4M\cite{sakkos09hoc}
and on the other hand, methods without modifying potentials with the reference methods V82 \cite{McL95a}, the fourth order complex triple-jump scheme (\ref{eq:TripleJump}), referenced as Yoshida 4 and a 6th-order complex coefficient method by Chambers \cite{chambers03siw}.
{We remark that all relevant methods in the cited papers have been tested and the most efficient ones for this problem are included in the plots.}

\subsubsection{P\"oschl-Teller potential} We have chosen
the well-known  one-dimensional P\"oschl-Teller potential for the
availability of analytic solutions of the eigenstates 
\beq\label{eq:poeschl-teller}
    H = -\frac12\frac{\partial^2}{\partial x^2} - \frac{\lambda(\lambda+1)}{2} \left(\sech(x)^2-1\right),
 \eeq
with $\lambda(\lambda+1)=10$.  
The results of our computation are shown in
Figure~\ref{fig:unperturbed}. The higher order of the complex
coefficient methods outweighs their extra cost starting from
moderate accuracy. The optimizations of the error terms can be
clearly appreciated in the comparison with the 4th order
triple-jump (\ref{eq:TripleJump}). When we consider methods with
modified potentials, we observe that the new methods show only
slight improvements with respect to the method OMF-4M since both
parts of the splitting $T$ and $V$ are of comparable size. As the
desired precision is increased, 
 the new sixth order methods dominate in efficiency.

\setlength{\figurewidth}{.41\textwidth}
\setlength{\figureheight}{.66\figurewidth}
    
\begin{figure*}
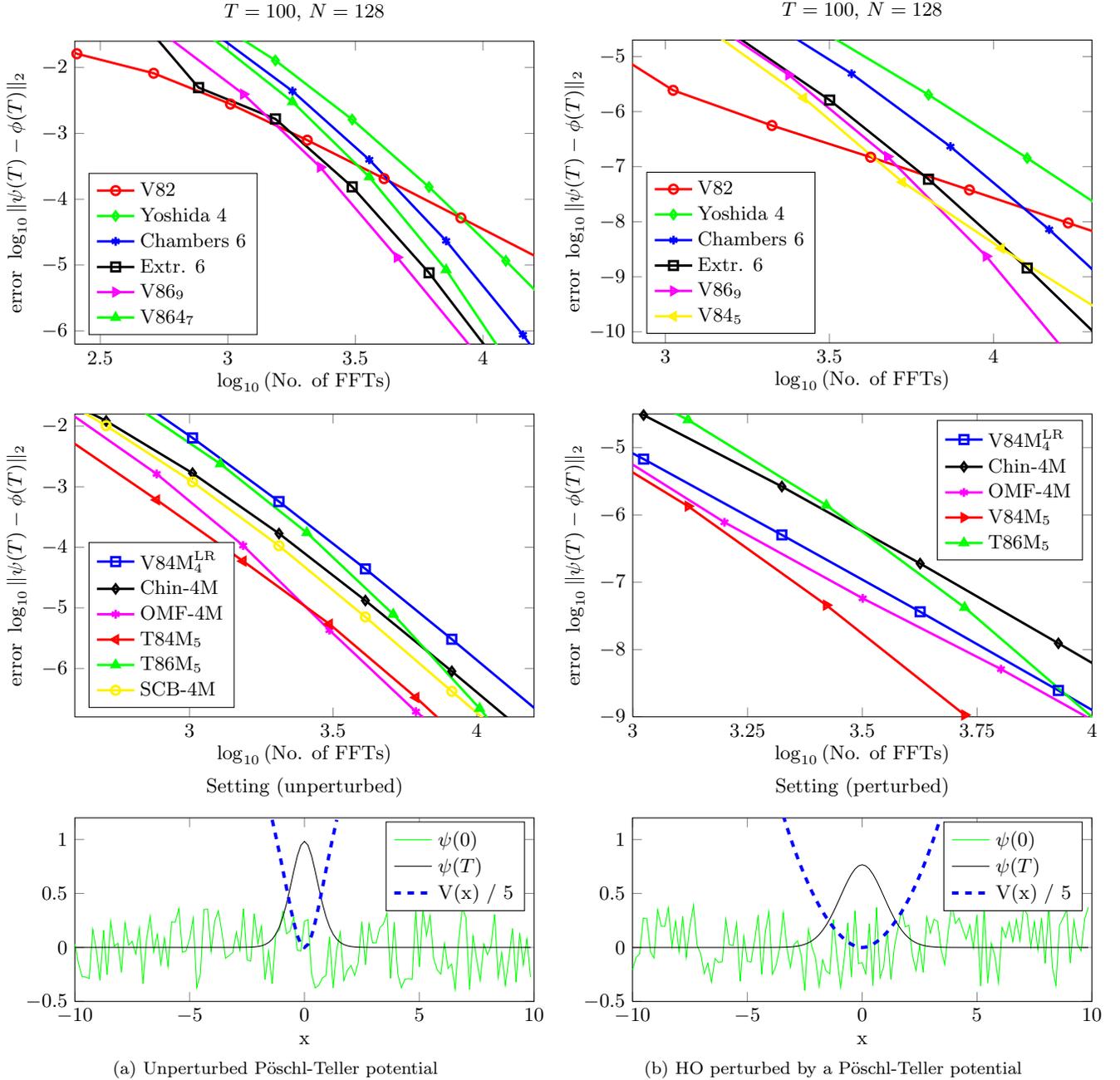
\centering
    \subfloat[Unperturbed P\"oschl-Teller potential]{\label{fig:unperturbed}%
{\definecolor{mycolor1}{rgb}{1,0,1}%
\begin{tikzpicture}
\begin{axis}[%
name=plot1,
width=\figurewidth,height=\figureheight,
unbounded coords=jump,
scale only axis,
xmin=2.4, xmax=4.2,
xlabel={$\log_{10}\text{(No. of FFTs)}$\rule[-3mm]{0mm}{3mm}},
xlabel style={yshift=2mm},
ymin=-6.2, ymax=-1.6,
ylabel={$\text{error }\log_{10}\|\psi(T)-\phi(T)\|_2$},
ylabel style={yshift=-3mm},
title={$T= 100$, $N= 128$},
legend style={at={(0.03,0.03)},anchor=south west,draw=black,fill=white,legend cell align=left},
xtick={0,0.5,1,1.5,2,2.5,3,3.5,4,4.5,5,5.5,6},
ytick={-13,-12,-11,-10,-9,-8,-7,-6,-5,-4,-3,-2,-1},
]
\input{plots_const_unperturbed1.tikz}
\end{axis}

\begin{axis}[%
name=plot2,
width=\figurewidth, height=\figureheight,
unbounded coords=jump,
scale only axis,
xmin=2.6, xmax=4.2,
xlabel={$\log_{10}\text{(No. of FFTs)}$},
xlabel style={yshift=2mm},
ylabel style={yshift=-3mm},
ymin=-6.8, ymax=-1.8,
ylabel={$\text{error }\log_{10}\|\psi(T)-\phi(T)\|_2$},
at=(plot1.below south west), anchor=above north west,
legend style={at={(0.03,0.03)},anchor=south west,draw=black,fill=white,legend cell align=left},
xtick={0,0.5,1,1.5,2,2.5,3,3.5,4,4.5,5},
ytick={-13,-12,-11,-10,-9,-8,-7,-6,-5,-4,-3,-2,-1},
]
\input{plots_const_unperturbed2.tikz}

\end{axis}

\begin{axis}[%
width=\figurewidth,
height=.4\figurewidth,
scale only axis, xmin=-10,xmax=10,
xlabel={x}, xlabel style={yshift=1mm},
ymin=-0.5,ymax=1.2,
at=(plot2.below south west),anchor=above north west,
title={Setting (unperturbed)},
legend style={draw=black,fill=white,legend cell align=left}
]
\addplot [
color=green,
solid
]
table[row sep=crcr]{
-10 0.0329373986487521\\
-9.84375 -0.146459139525481\\
-9.6875 -0.269691576085853\\
-9.53125 -0.275003887598898\\
-9.375 -0.287405337564394\\
-9.21875 0.166089144660096\\
-9.0625 -0.0278147955008901\\
-8.90625 -0.306151116773819\\
-8.75 0.159024314497519\\
-8.59375 -0.253336697887454\\
-8.4375 0.240396675770248\\
-8.28125 0.0110482510912495\\
-8.125 0.038345884542549\\
-7.96875 -0.231282617974489\\
-7.8125 0.225292758189547\\
-7.65625 0.0209703293127552\\
-7.5 0.0562251405451114\\
-7.34375 -0.0617158787913845\\
-7.1875 0.175076819580893\\
-7.03125 -0.337922019538481\\
-6.875 0.0750967628099589\\
-6.71875 0.286556805711229\\
-6.5625 -0.0405312761956822\\
-6.40625 0.120786805816685\\
-6.25 -0.155579107697148\\
-6.09375 0.0850540729122744\\
-5.9375 -0.175030413592702\\
-5.78125 0.237146395720145\\
-5.625 0.234460175451217\\
-5.46875 0.35947466698509\\
-5.3125 -0.0440642954717268\\
-5.15625 -0.034133414654934\\
-5 0.0790195274919656\\
-4.84375 0.271234135689711\\
-4.6875 -0.37112182667863\\
-4.53125 -0.24757036782803\\
-4.375 0.351168350798316\\
-4.21875 0.354594374914046\\
-4.0625 -0.0372195722648619\\
-3.90625 0.246068378058379\\
-3.75 0.339519168468027\\
-3.59375 0.136730406893578\\
-3.4375 -0.101067242676091\\
-3.28125 -0.0746553354136824\\
-3.125 -0.0484296071456034\\
-2.96875 0.141426368939748\\
-2.8125 -0.0276518257826484\\
-2.65625 0.358823354463861\\
-2.5 -0.11502730885578\\
-2.34375 -0.127437989612105\\
-2.1875 0.313379790086996\\
-2.03125 0.0359677073598597\\
-1.875 0.19734368157949\\
-1.71875 -0.296963296917258\\
-1.5625 -0.0370213235763228\\
-1.40625 -0.336646582084206\\
-1.25 0.129307434849829\\
-1.09375 0.161219373683976\\
-0.9375 0.331659237253091\\
-0.78125 0.126719787635769\\
-0.625 0.150495207063705\\
-0.46875 0.28002277954233\\
-0.3125 -0.0254103729102257\\
-0.15625 -0.0328700561023245\\
0 0.242325202563447\\
0.15625 0.257099672374732\\
0.3125 -0.245064056870439\\
0.46875 -0.375497430575434\\
0.625 -0.350844283986951\\
0.78125 -0.28266750793171\\
0.9375 -0.260119019096702\\
1.09375 0.0996241259605968\\
1.25 -0.372455939075509\\
1.40625 -0.0219029654051167\\
1.5625 0.141262949832731\\
1.71875 -0.304949572413816\\
1.875 -0.20894311709478\\
2.03125 -0.1669728409604\\
2.1875 -0.259060151815075\\
2.34375 -0.139561793908938\\
2.5 0.238371667960935\\
2.65625 -0.158618556736357\\
2.8125 0.218197143103763\\
2.96875 0.0417756526415506\\
3.125 0.0432963123840009\\
3.28125 0.18259028379189\\
3.4375 0.216613144330898\\
3.59375 0.31732775900628\\
3.75 -0.286442396822363\\
3.90625 0.232803819815986\\
4.0625 -0.245873692981198\\
4.21875 -0.372884079504709\\
4.375 -0.294976278545031\\
4.53125 -0.289950538574533\\
4.6875 -0.29424498479105\\
4.84375 0.344570318988665\\
5 -0.179520813738647\\
5.15625 0.350443884064401\\
5.3125 0.109365861775158\\
5.46875 0.294911374403112\\
5.625 -0.105230525423317\\
5.78125 -0.208833885420735\\
5.9375 -0.247532292578281\\
6.09375 0.0361379158596423\\
6.25 -0.193860316413427\\
6.40625 -0.153718359631354\\
6.5625 -0.383512490058388\\
6.71875 0.0692593271450419\\
6.875 0.366180996482166\\
7.03125 0.276954422679085\\
7.1875 -0.389535563612759\\
7.34375 0.106100364688509\\
7.5 -0.111391171926783\\
7.65625 -0.305526960229834\\
7.8125 0.0323132210727561\\
7.96875 -0.0661644174297029\\
8.125 0.0135359038459025\\
8.28125 0.305671428447308\\
8.4375 -0.277553235239876\\
8.59375 -0.0517128661652328\\
8.75 -0.349083790737077\\
8.90625 -0.094181709142572\\
9.0625 0.176037207846704\\
9.21875 -0.320520982345121\\
9.375 0.132355088848168\\
9.53125 -0.16117805285391\\
9.6875 0.0780265110393289\\
9.84375 -0.275608895769459\\
};
\addlegendentry{$\psi(0)$};

\addplot [
color=black,
solid
]
table[row sep=crcr]{
-10 7.84803006254201e-012\\
-9.84375 -8.01301059928027e-012\\
-9.6875 8.09870696992595e-012\\
-9.53125 -8.1939956556865e-012\\
-9.375 8.28042135269311e-012\\
-9.21875 -8.37167235347374e-012\\
-9.0625 8.45685691145214e-012\\
-8.90625 -8.54384538062311e-012\\
-8.75 8.63051042390382e-012\\
-8.59375 -8.71766917311757e-012\\
-8.4375 8.80248265827895e-012\\
-8.28125 -8.89002641920537e-012\\
-8.125 8.97641629802993e-012\\
-7.96875 -9.0617334297662e-012\\
-7.8125 9.14840522194426e-012\\
-7.65625 -9.23265074082153e-012\\
-7.5 9.32021736802718e-012\\
-7.34375 -9.40159186677704e-012\\
-7.1875 9.49115528863099e-012\\
-7.03125 -9.57025585814744e-012\\
-6.875 9.65813384923519e-012\\
-6.71875 -9.7363470236869e-012\\
-6.5625 9.82490253408974e-012\\
-6.40625 -9.90049860065678e-012\\
-6.25 9.99421726881165e-012\\
-6.09375 -1.00304777786532e-011\\
-5.9375 1.02808748068911e-011\\
-5.78125 -9.69080865906213e-012\\
-5.625 1.23326787644573e-011\\
-5.46875 -2.87111923560689e-012\\
-5.3125 3.72680303307849e-011\\
-5.15625 8.21901600811798e-011\\
-5 3.20974041288383e-010\\
-4.84375 9.95512365881313e-010\\
-4.6875 3.17010816753155e-009\\
-4.53125 9.59916892119373e-009\\
-4.375 2.8332210958256e-008\\
-4.21875 8.08704942811245e-008\\
-4.0625 2.23887123304286e-007\\
-3.90625 6.00715241948094e-007\\
-3.75 1.56297099287254e-006\\
-3.59375 3.94383556119287e-006\\
-3.4375 9.65380992446948e-006\\
-3.28125 2.29294164584918e-005\\
-3.125 5.28592026911339e-005\\
-2.96875 0.00011830378662813\\
-2.8125 0.000257126225291992\\
-2.65625 0.00054284721890621\\
-2.5 0.00111351274127357\\
-2.34375 0.00221962769767089\\
-2.1875 0.00430008870370785\\
-2.03125 0.00809608891763016\\
-1.875 0.01481075296311\\
-1.71875 0.0263126149716105\\
-1.5625 0.0453571665129535\\
-1.40625 0.0757534061181663\\
-1.25 0.122324287600671\\
-1.09375 0.190409017237213\\
-0.9375 0.284582734003898\\
-0.78125 0.406356828749265\\
-0.625 0.551076776827864\\
-0.46875 0.705168170890212\\
-0.3125 0.845919631655522\\
-0.15625 0.945967579660391\\
0 0.98233128461437\\
0.15625 0.94596757966035\\
0.3125 0.845919631655461\\
0.46875 0.705168170890155\\
0.625 0.551076776827823\\
0.78125 0.406356828749242\\
0.9375 0.284582734003877\\
1.09375 0.190409017237201\\
1.25 0.122324287600666\\
1.40625 0.075753406118167\\
1.5625 0.0453571665129563\\
1.71875 0.0263126149716242\\
1.875 0.0148107529631291\\
2.03125 0.00809608891764533\\
2.1875 0.00430008870371257\\
2.34375 0.00221962769767078\\
2.5 0.0011135127412754\\
2.65625 0.000542847218909931\\
2.8125 0.000257126225292832\\
2.96875 0.000118303786635305\\
3.125 5.28592027042417e-005\\
3.28125 2.29294164690064e-005\\
3.4375 9.65380993506019e-006\\
3.59375 3.94383556978464e-006\\
3.75 1.56297099704849e-006\\
3.90625 6.00715245022133e-007\\
4.0625 2.23887125349517e-007\\
4.21875 8.08704950181671e-008\\
4.375 2.8332211833358e-008\\
4.53125 9.59916710363553e-009\\
4.6875 3.17010186966691e-009\\
4.84375 9.95512666265929e-010\\
5 3.20975321316453e-010\\
5.15625 8.21915431786934e-011\\
5.3125 3.7267648417725e-011\\
5.46875 -2.86893790962632e-012\\
5.625 1.23348860708867e-011\\
5.78125 -9.6901275415701e-012\\
5.9375 1.02777549895776e-011\\
6.09375 -1.00304110118421e-011\\
6.25 9.99502747780999e-012\\
6.40625 -9.89832323543485e-012\\
6.5625 9.82605782813317e-012\\
6.71875 -9.73503722028794e-012\\
6.875 9.65959792855022e-012\\
7.03125 -9.57045733889176e-012\\
7.1875 9.48781082090235e-012\\
7.34375 -9.40354361148342e-012\\
7.5 9.32028750122788e-012\\
7.65625 -9.23121382800987e-012\\
7.8125 9.14887103881463e-012\\
7.96875 -9.06089886050765e-012\\
8.125 8.97736284533571e-012\\
8.28125 -8.89001250430285e-012\\
8.4375 8.80328850557603e-012\\
8.59375 -8.71654273108725e-012\\
8.75 8.63069412014133e-012\\
8.90625 -8.54400475262006e-012\\
9.0625 8.45552622475248e-012\\
9.21875 -8.37164485551749e-012\\
9.375 8.2792015889767e-012\\
9.53125 -8.19554738884057e-012\\
9.6875 8.098900885659e-012\\
9.84375 -8.01204148845018e-012\\
};
\addlegendentry{$\psi(T)$};

\addplot [
color=blue,
dashed,
line width=1.5pt
]
table[row sep=crcr]{
-1.40625 1.18173610044121\\
-1.25 1.03208513381957\\
-1.09375 0.876449307675426\\
-0.9375 0.714642245894227\\
-0.78125 0.549032698039512\\
-0.625 0.385705852030012\\
-0.46875 0.235079346353686\\
-0.3125 0.111164430232318\\
-0.15625 0.0289049390011292\\
0 0\\
0.15625 0.0289049390011292\\
0.3125 0.111164430232318\\
0.46875 0.235079346353686\\
0.625 0.385705852030012\\
0.78125 0.549032698039512\\
0.9375 0.714642245894227\\
1.09375 0.876449307675426\\
1.25 1.03208513381957\\
1.40625 1.18173610044121\\
};
\addlegendentry{V(x) / 5};

\end{axis}

\end{tikzpicture}%
}
    }
        \subfloat[HO perturbed by a P\"oschl-Teller potential]{\label{fig:perturbed}%
        \hspace{0cm}
{\input{const_perturbed.tikz}}
	  }
    \caption{\label{fig:fixedstep} (color online)
    In the first row, efficiency curves (error vs. number of FFTs) for methods without force evaluations are presented, with the new methods (triangles) performing best for high accuracies.
    The middle rows depicts methods based on modified potentials.
    In the right column, T86M$_5$ intersects with V84M$_5$ at precision $10^{-13}$, whereas it already improves on T84M$_5$ at $10^{-9}$ for the left column. SCB-4M overlaps with Chin-4M and has thus been omitted in the plot.
        In the bottom row, the random initial conditions (green), the ground
states (black) and the potentials (dashed blue), scaled by $1/5$
to fit the axis, are shown. }
\end{figure*}

\subsubsection{Perturbed harmonic oscillator} To
illustrate the benefits of methods designed for near integrable
systems, we use the Hamiltonian
$$
    H = -\frac12\frac{\partial^2}{\partial x^2} + \frac12 \omega^2 x^2 + \varepsilon V_{\varepsilon}(x),
$$
and split it in a large part $H_{\mathrm{HO}} =
-\frac12\frac{\partial^2}{\partial x^2} + \frac12\omega^2 x^2$ and a
small part $\varepsilon V_{\varepsilon}(x)$. The trap frequency is
set to $\omega=1$ and the perturbation $\varepsilon V_\varepsilon$
is given by the P\"oschl-Teller potential in
\eqref{eq:poeschl-teller}, with $\lambda(\lambda+1)=2/5$.
 The harmonic
part $H_{\mathrm{HO}}$ can be solved exactly via an exact
splitting using Fourier transforms, cf. {\cite{chin05foa}},
where it is shown that
$$
    e^{-i\delta H_{\mathrm{HO}}} \equiv
                        e^{-i \frac{\omega}{2}\tan(\frac{\delta\omega}{2})\, x^2}\,
                        e^{-i\frac1{2\omega}\sin(\delta\omega)\, p^2}\,
                        e^{-i\frac{\omega}{2}\tan(\frac{\delta\omega}{2})\, x^2},
$$
for  $|\delta\omega|<\pi$ and
$p^2\equiv-\frac{\partial^2}{\partial x^2}$.

From the computational point of view, it is suggested \cite{bader11fmf} to
  consider the VTV split instead of the TVT split because it can be
  concatenated with the perturbation which only depends on the
  coordinates and no additional FFTs are necessary, i.e. 
$$
 \ldots e^{-ib_{j+1} \tau \varepsilon V}
   e^{-ia_j \tau H_{\mathrm{HO}}}
  e^{-ib_j \tau \varepsilon V}  \ldots
$$
In \cite{chin05foa} this decomposition is
  generalized to the two-dimensional problem
  $H=\frac12(p_x^2+p_y^2)+\frac12(w_1^2x^2+w_2^2y^2)-\Omega(xp_y-yp_x)$
  and in \cite{bader11fmf} to the non homogeneous and possibly time-dependent
  one-dimensional problem $H=\frac12p^2+\frac12 w(t)x^2 + f(t)p+ g(t)x$.

After the substitution $\delta=-ih$, we have
$$
    e^{-hH_{\mathrm{HO}}} \equiv
e^{-\frac{\omega}{2}\tanh(\frac{h\omega}{2})\, x^2}\,
e^{-\frac{1}{2\omega}\sinh(h\omega)\, p^2}\,
e^{-\frac{\omega}{2}\tanh(\frac{h\omega}{2})\, x^2},
$$
for $|\Im(h)\omega|<\pi$ and $\Re(h)>0$ (for numerical stability) and the perturbation part is easily propagated after discretization by the exponential of a diagonal matrix.
In this setting, the higher order in the small parameter is amplified and the efficiency plots in Figure~\ref{fig:perturbed} indicate that the new methods outperform the existing ones when high precision is sought and overall when modified potentials are allowed.
We observe in both examples that, when modified potentials can be computed without exceedingly large computational cost, they should be used.


\structurecomment{Remark on mesh size} Further numerical
experiments show that the efficiency curves are independent of the
mesh size, i.e., the norm of $T$, and the cost only increases as
$N\log(N)$ as expected. The reason for this can be understood by
following the evolution of the state vector along the iterations
of the algorithm. Whereas in the beginning one has a non-smooth
configuration $u_0$, after a few steps the vector $u_i$ is close
to an eigenstate and thus smoothened.

It is important to remark that the methods proposed in this work
can be implemented in an algorithm which uses variable step,
variable order, variable mesh size and variable simple-double
precision. The best implementation can depend on the class of
problems to be solved. For illustration, we present an
implementation with variable time steps.

\subsection{Variable step method}\label{subsec:varstep}

The previous examples show that for low accuracies and large time
steps, the (8,2) method (with real coefficients) performs best.
However, if we allow for variable time steps, as proposed in
\cite{aichinger05afc,lehtovaara07sot}, the computational cost is
drastically reduced. We propose an improved time-stepping
algorithm that is based on two different estimators for the
eigenvalue.

Recall that fixing the time-step and iterating to convergence will
yield an eigenvector with the error being of the order of the
method $\mathcal{O}(h^{p})$ since we are computing exactly the
spectrum of a perturbed Hamiltonian. Assume now that we are close
to convergence, i.e, one has obtained an eigenvector $u_n = v_0 +
\mathcal{O}\left(h^{p}\right)$ and we consider the decomposition
in the basis of exact eigenvectors $v_i$ of $H$,
$$
  u_n  = \sum_{i=0}^{N-1} d_i v_i, \quad \text{where}\quad \sum_{i=0}^{N-1} |d_i|^2 =1.
$$

It is clear that $d_i=\cO(h^{p}), \ i>1$ and due to the
normalization $d_0=1+\cO(h^{2p)})$. Then, an energy estimation is
given by
\begin{equation*}
  E_{h,1}\equiv u_n^T H u_n
  = E_0 + \cO(h^{2p}).
\end{equation*}
Alternatively, the energy can be estimated by the loss of norm in
each time step,
\[
 \bar u_{n+1} 
                            = e^{-hH}u_n + \cO{(h^{p+1})}
 = e^{-hE_{0}}v_0 + \cO{(h^{p+1})},
\]
and then
\begin{align*}
    E_{h,2} &\equiv \frac{\log\left(\left\|\bar u_{n+1}\right\|\right)}{h} =  E_{0} + c h^{p} + \cO\left(h^{p+1}\right).
\end{align*}
Combining both expressions yields an error estimate for the
energy,
\begin{equation*}
    \Delta E_h \equiv E_{h,2} - E_{h,1} = c h^{p} + \cO(h^{p+1}).
\end{equation*}
The convergence in energy is measured by comparison with the
previous time step,
$$
    \boldsymbol{\delta}
       E_h^n \equiv E_{h,1}^{n}-E_{h,1}^{n-1} = d h^{2p}+\cO(h^{2p+1}).
$$
The time stepper then works as follows: starting from a large
step size, the time step is decreased by a factor $1/2$ whenever
the actual reduction in energy of the iteration
$\boldsymbol{\delta}E$ falls below the the maximally reachable
precision $\Delta E$, i.e., $|\boldsymbol{\delta} E| < (\Delta E)^2$
and the iteration is terminated once the error estimate $\Delta E$
has reached a given tolerance.

\begin{figure*}
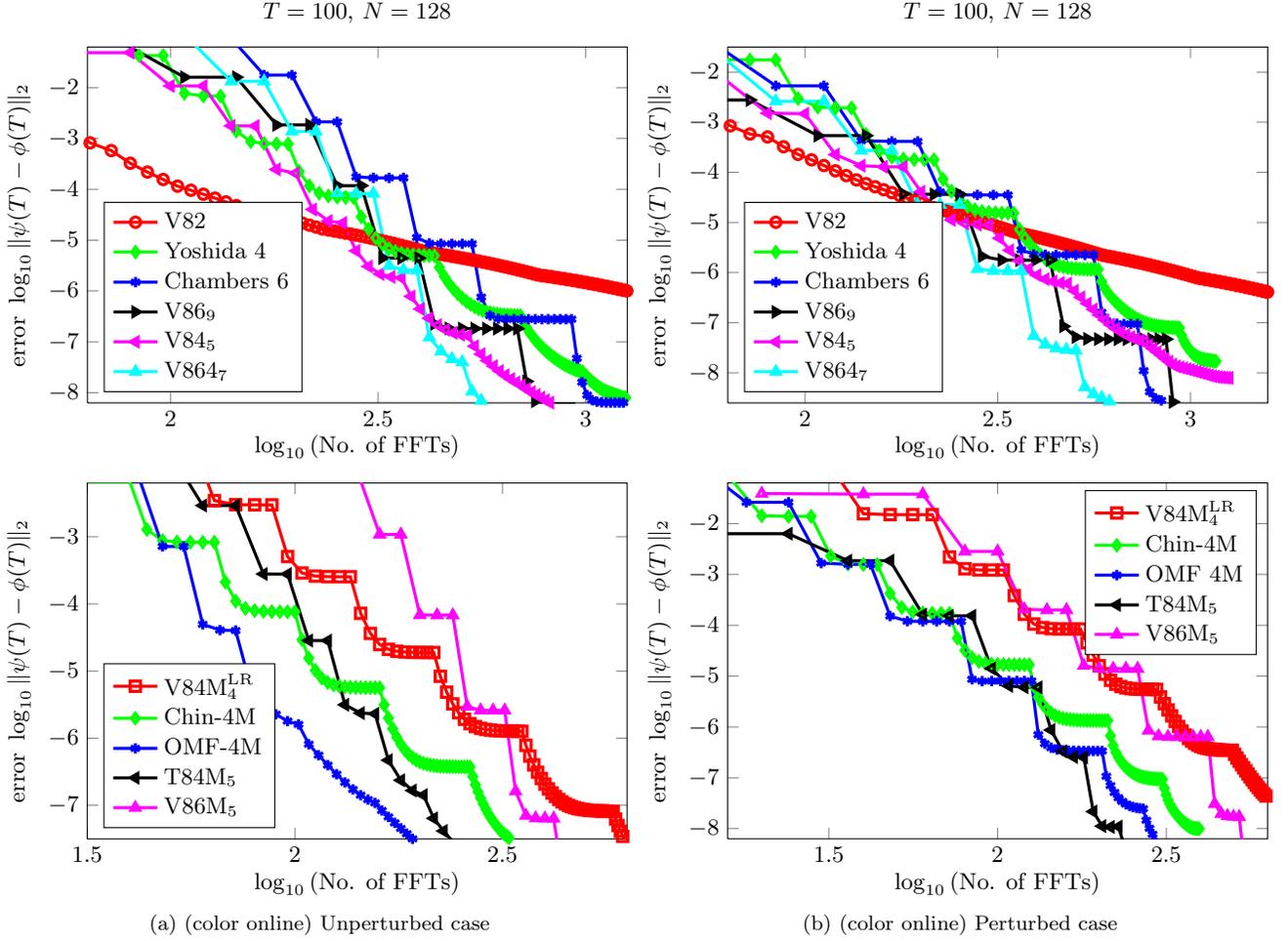
\centering
    \subfloat[\label{fig:varstep_unperturbed}%
        (color online) Unperturbed case]{
			{\input{var_unperturbed.tikz}}
    }
    \subfloat[\label{fig:varstep_perturbed}%
        (color online) Perturbed case]{%
       {\input{var_perturbed.tikz}} 
    }
    \caption{\label{fig:varstep}Evolution of precision in the $\ell_2$ norm of the position vector with the variable time step algorithm described in Sec. \ref{subsec:varstep}. As in Fig.~\ref{fig:fixedstep}, the top row gives the results for standard methods whereas the bottom rows shows methods with modifying potentials.
    }
\end{figure*}

For the numerical
experiments, we use the same configurations as for constant time
step but terminate the algorithm when convergence in
energy is reached at $\Delta E<10^{-10}$. The iterations are
initialized with random normalized data and a time step of
$\tau=10$.
The results are displayed in Figure~\ref{fig:varstep_unperturbed} for the P\"oschl-Teller potential and in Figure~\ref{fig:varstep_perturbed} for the
perturbed harmonic oscillator with the same parameters as in the fixed-step size experiments.
The error is measured as the $\ell_2$ norm of the difference between the current value of the algorithm $\psi(t)$ and the exact ground state $\phi(T)$ as in the previous experiments, $\text{error} = \|\psi(t)-\phi(T)\|$.


\structurecomment{discussion of result} As expected, it is
apparent that lower order methods show better smoothing behavior
for the first steps, when the wave function is still rough (recall
that the algorithm is initialized by a worst-case wave function).
For higher precisions, the new methods clearly outperform the
existing ones, with the sole exception of the unperturbed setting
with modified potentials, where the globally optimized OMF-4M
method can hardly be improved unless extremely high precision is
sought and the 6th order methods of Table~\ref{tab.2a} and
\ref{tab.2b} become favorable (not shown).
Finally, if one is interested in very high
  accuracies, high order extrapolation methods \cite{chin09aoi,blanes99eos} can be used for the last part of the time integration.
 

\structurecomment{Variable order - or put it with conclusions?}
The results indicate that for low precision, i.e., for the first
iterations, a lower order method should be used and then, after a
certain precision is reached, e.g., when the higher order methods exhibit their superiority
the algorithm should change to the optimal method, either V864$_7$
or V86M$_5$ until convergence.
Further preliminary experiments on this adaptive order strategy have shown that there is plenty of room for optimization, e.g., by changing the initial step-size, adjusting the step-size by a different factor or by modifying the control criterion. Each of which has certain advantages and disadvantages, depending on the initial conditions and the range of precision.

\structurecomment{On excited states} For excited states, one
expects an even better performance of the new methods since
several states have to be computed to high precision in order to
avoid error accumulation and the gains of the new methods are thus
amplified. We have  confirmed this conjecture by numerical
experiments.
The results thereof are omitted in the manuscript since they do not contribute insight beyond the presented experiments: they are qualitatively identical.
%
%
%

\section{Conclusions and outlook}
We have studied the Schr\"odinger eigenvalue problem by the imaginary time propagation method and proposed splitting schemes with
positive real coefficients using modified potentials as well as with 
 complex coefficients that can overcome the order barrier for parabolic problems since the coefficients have only positive real parts.
The obtained sixth order methods are clearly superior to any classical ones for high precisions.
On the other hand, when the gradient of the potential can be cheaply evaluated, the high order methods with complex coefficients are efficient only at very high accuracies due to the double cost caused by complex arithmetic.

We have proposed different high order methods to reach highly accurate results.
An efficient implementation should take into account, for example, a preliminary time integration on a coarse mesh using simple precision arithmetic in order to get, as fast as possible, a smooth and relatively accurate solution from a random initial guess, and next consider a refined mesh using arithmetic in double precision.
For simple precision arithmetic and low accuracies, it suffices to consider only low order methods, and when higher accuracies are desired we turn to double precision, variable time step and variable order methods.
The best algorithm could depend on the class of problems to solve.

It is also important to remark that the form of the exponent allows that the techniques presented in this work can also be transferred to other areas whenever splitting is appropriate and the integration has to be performed forward in time, e.g., statistical mechanics of quantum systems, where one has to compute the Boltzman operator $\exp(- \beta H)$, with $\beta =(kT)^{-1}$ or quantum Monte-Carlo simulations \cite{bandrauk06cis}.

 Finally, we would like to mention that real time integration
 with complex coefficients is under investigation. To compute
 $\e^{-ia_itT}$ requires complex FFTs, and this is irrespective of
 the coefficients $a_i$ being real or complex. However, the constraint
 $\Im(a_i)\leq 0$ and the consistency condition, $\sum_ia_i=1$, necessarily
 requires $a_i\in\mathbb{R}$, while $b_i$ can be complex. A
large number of new methods have been explored, but the
superiority is not yet clear since there exist highly efficient
methods with real coefficients for perturbed problems
\cite{blanes13nfo,celesmech2000} and using modified potentials
\cite{blanes01hor}. 

\acknowledgments{
We wish to acknowledge Ander Murua and Joseba Makazaga for providing the methods T86$_9$ and V86$_9$. This work has
been partially supported by Ministerio de Ciencia e Innovaci\'on (Spain) under project MTM2010-18246-C03 and by a 
grant from the Qatar National Research Fund \#NPRP  NPRP 5-674-1-114.
P.B. also acknowledges the support through the FPU fellowship AP2009-1892.
}


\end{document}